\documentclass{amsart}

\usepackage{amssymb}
\usepackage{graphicx}
\usepackage{amsmath}
\usepackage{amsfonts}
\usepackage[usenames]{color}
\usepackage{latexsym}
\usepackage{array}
\usepackage[all]{xy}
\usepackage{verbatim}
\usepackage{comment}
\input epsf

\newtheorem{theorem}{Theorem}[section]

\newtheorem{corollary}{Corollary}[section]

\newtheorem{lemma}{Lemma}[section]

\theoremstyle{remark}

\newtheorem*{question}{Question}

\newtheorem{definition}{Definition}[section]
\newtheorem*{ConGor}{\bf{The Conway-Gordon Theorem}}

\makeatletter
 
 \@addtoreset{equation}{section}
\makeatother

\def\Z{{\mathbb Z}}

\def\R{{\mathbb R}}

\def\TSG{{\mathrm{TSG_+}}}
\def\Aut{{\mathrm{Aut}}}

\def\so{{\mathrm{SO}}}

\def\TSG{{\mathrm{TSG}_+}}

\def\Aut{{\mathrm{Aut}}}

\newcommand{\x}{\times}

\newcommand{\dfn}[1]{\textit{#1}}

\newcommand{\ty}{\nabla\mathrm{Y}}
\newcommand{\yt}{\mathrm{Y}\nabla}

\numberwithin{equation}{section}

\begin{document}
\title{Recent Developments in Spatial Graph Theory}


\author[E. Flapan]{Erica Flapan}
\address{Department of Mathematics, Pomona College, Claremont, CA 91711, USA}
\email{eflapan@pomona.edu}

\author[T. Mattman]{Thomas W.\ Mattman}
\address{Department of Mathematics and Statistics,
California State University, Chico,
Chico, CA 95929-0525}
\email{TMattman@CSUChico.edu}

\author[B. Mellor]{Blake Mellor}
\address{Department of Mathematics, Loyola Marymount University, Los Angeles, CA 90045, USA}
\email{blake.mellor@lmu.edu}

\author[R. Naimi]{Ramin Naimi}
\address{Occidental College, Los Angeles, CA 90041}
\email{rnaimi@oxy.edu}

\author[R. Nikkuni]{Ryo Nikkuni}
\address{Department of Mathematics, School of Arts and Sciences, Tokyo Woman's Christian University, 2-6-1 Zempukuji, Suginami-ku, Tokyo 167-8585, Japan}
\email{nick@lab.twcu.ac.jp}


\date \today

\begin{abstract} This article presents a survey of some recent results in the theory of spatial graphs.  In particular, we highlight results related to intrinsic knotting and linking and results about symmetries of spatial graphs.  In both cases we consider spatial graphs in $S^3$ as well as in other $3$-manifolds.

\end{abstract}

\maketitle

\section{Introduction}\label{intro}   Spatial graph theory is the study of graphs embedded in $S^3$.  Much of the current work in this area has its roots in John Conway and Cameron Gordon's \cite{CG} result from 1983 that every embedding of the complete graph $K_6$ in $S^3$ contains a non-split link and every embedding of $K_7$ in $S^3$ contains a non-trivial knot (note the result about $K_6$ was independently obtained by Horst Sachs \cite{Sa1,Sa2}) .  Because these properties are independent of the way that $K_6$ or $K_7$ are embedded in $S^3$, we say the properties are {\it intrinsic} to the graph.  Thus $K_6$ is said to be {\it intrinsically linked} and $K_7$ is said to be {\it intrinsically knotted}.   Conway and Gordon's theorem has motivated the study of a variety of intrinsic properties of graphs which will be discussed in Sections 2--6.

Independent of Conway and Gordon's ground breaking results on intrinsic knotting and linking, spatial graph theory also has roots in the study of symmetries of non-rigid molecules.  Chemists have long used the symmetries of a molecule to characterize it and predict its properties. For small molecules, it is enough to consider rotations and reflections.  Increasingly however, chemists are working with large molecules ranging from synthetic structures in the form of a knot, link, or M\"{o}bius ladder to polymers such as proteins and DNA. Because of their size, these molecules can be relatively flexible, and hence their symmetries cannot always be seen as rigid motions. To characterize the symmetries of such complex molecules, Jon Simon introduced the {\em topological symmetry group} \cite{si}. Although the motivation for studying this group comes from molecular symmetries, the topological symmetry group provides useful information about any spatial graph in $S^3$.  In Section~\ref{symmetries}, we present results about topological symmetry groups and other measures of symmetry of spatial graphs in $S^3$.

In addition to studying spatial graphs in $S^3$, it is natural to consider embeddings of graphs in any $3$-manifold.  In Section~\ref{3man}, we discuss intrinsic knotting and linking of spatial graphs in arbitrary $3$-manifolds as well as symmetry and asymmetry of spatial graphs in $3$-manifolds.

\medskip

\section{Intrinsic linking and knotting}\label{IKIL}
For completeness, we begin with several definitions.

\begin{definition}
Let $G$ be an abstract graph and let $f$ be an embedding of $G$ in $S^3$.  Then we say the image $f(G)$ is a {\it spatial graph} and $f$ is a {\it spatial embedding}.  \end{definition}

\begin{definition}  If the image of every embedding of $G$ in $S^3$ contains a non-split link then we say $G$ is {\it intrinsically linked}, and if the image of every embedding of $G$ in $S^3$ contains a non-trivial knot then we say $G$ is {\it intrinsically knotted}.  

\end{definition}

\begin{definition} A \dfn{minor} of graph is a graph obtained by contracting zero or more edges of a subgraph.  We say that a graph $G$ is \dfn{minor minimal} with respect
to a property ${\mathcal P}$ if $G$ has ${\mathcal P}$ but no proper minor of $G$ does. 
\end{definition}

Intuitively, we can think of a minor of a graph $G$ as a graph obtain by deleting and/or contracting edges of $G$.

As mentioned in the introduction, Conway and Gordon \cite{CG} proved that $K_6$ is intrinsically linked and $K_7$ is intrinsically knotted.  It is not hard to check that every minor of $K_6$ has a linkless embedding and every minor of $K_7$ has a knotless embedding.  Thus $K_6$ is  minor minimal with respect to being intrinsically linked and $K_7$ is minor minimal with respect to be being intrinsically knotted.  

\begin{figure}[h]
{\includegraphics{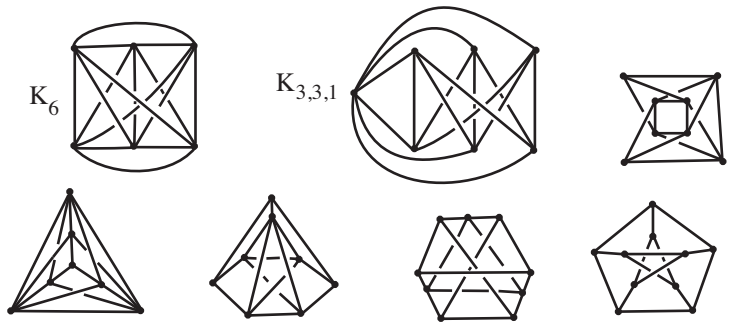}}
\caption{The Petersen family of graphs. }
\label{Petersen}
\end{figure}

At the same time that Conway and Gordon proved their theorem, Horst Sachs \cite{Sa1}, \cite{Sa2} independently proved that every graph in the Petersen family (illustrated in Figure~\ref{Petersen}) is minor minimal with respect to being intrinsically linked.  This includes the graphs $K_6$ and $K_{3,3,1}$ as well as five other graphs.  Subsequently, Robertson, Seymour, and Thomas~\cite{rst} obtained the surprising result that the graphs in the Peterson family are the only graphs which are minor minimal  with respect to being intrinsically linked.  Thus minor minimal intrinsically linked graphs are completely characterized.

  Robertson and Seymour's Graph Minor Theorem \cite{rs} implies that the set of graphs which are minor minimal with respect to being intrinsically knotted is finite.   Since any graph which contains an intrinsically knotted graph as a minor is itself intrinsically knotted, finding this finite list would enable us to determine whether or not any given graph is intrinsically knotted.  However, as of now, there is no known list of all minor minimal intrinsically knotted graphs.  Below we present a survey of graphs which are known to be minor minimal intrinsically knotted.
  
We begin by observing that the graphs in the Petersen family are related by two operations.   A {\it $\triangle Y$ move} is an operation to obtain a new graph $G$ from a graph $H$ by removing all edges of a $3$-cycle of $H$, and adding a new vertex and connecting it to each of the vertices of the cycle. A {\it $Y \triangle$ move} is the reverse of this operation.  Observe that neither of these moves changes the total number of edges in a graph.  It turns out that all of the graphs in the Petersen family can be obtained from either $K_6$ or $K_{3,3,1}$ by repeatedly applying  $\triangle Y$ and $Y \triangle$ moves.  Furthermore, the Petersen family is closed under $\triangle Y$ and $Y\triangle$ moves.  

It was shown in \cite{MRS88} that if we start with an intrinsically knotted graph and apply $\triangle Y$ moves we obtain other intrinsically knotted graphs. On the other hand, $Y \triangle$ moves do not necessarily preserve intrinsic knotting \cite{FN08}.

\begin{definition}We call the set of all graphs obtained from a graph $H$ by a finite sequence of $\triangle Y$ and $Y \triangle$ moves the {\it $H$-family} and denote it by ${\mathcal F}\left(H\right)$. The set of all graphs obtained from $H$ using only $\triangle Y$ moves is denoted by ${\mathcal F}_{\triangle}\left(H\right)$.  
\end{definition}

There are 14 graphs in ${\mathcal F}_{\triangle}\left(K_7\right)$, and Kohara and Suzuki~\cite{KS92} proved that all are minor minimal with respect to being intrinsically knotted.  There are six graphs in ${\mathcal F}\left(K_{7}\right)\setminus{\mathcal F}_\triangle\left(K_{7}\right)$.  However, it is known that none of these graphs is intrinsically knotted \cite{FN08}, \cite{GMN}, \cite{ HNTY10}.  Foisy \cite{F02} proved that the four-partite graph $K_{3,3,1,1}$ is intrinsically knotted, and  Kohara and Suzuki~\cite{KS92} proved that the 26 graphs in ${\mathcal F}_{\triangle}\left(K_{3,3,1,1}\right)$ are all minor minimal with respect to being intrinsically knotted.  In contrast with the graphs in ${\mathcal F}\left(K_{7}\right)\setminus{\mathcal F}_\triangle\left(K_{7}\right)$ which are not intrinsically knotted, there are 32 graphs in ${\mathcal F}\left(K_{3,3,1,1}\right)\setminus{\mathcal F}_\triangle\left(K_{3,3,1,1}\right)$ all of which turn out to be minor minimal intrinsically knotted \cite{GMN}.  Thus altogether there are 58 graphs in ${\mathcal F}\left(K_{3,3,1,1}\right)$ and all are minor minimal intrinsically knotted.

One approach to finding additional minor minimal intrinsically knotted graphs is to consider graphs with a small number of edges.  The 14 graphs in ${\mathcal F}_\triangle\left(K_7\right)$ 
all have exactly 21 edges.  In fact, the following theorem shows that these are the only intrinsically knotted graphs with 21 or fewer edges.

\begin{theorem}\cite{BM}, \cite{JKM}, \cite{LKLO}, \cite{Ma}
\label{thmminor minimal intrinsically knotted21}%
Every intrinsically knotted graph has at least 21 edges, and the graphs in ${\mathcal F}_\triangle\left(K_7\right)$  are the only intrinsically knotted graphs with exactly 21 edges.
\end{theorem}

\begin{figure}[h]
{\includegraphics{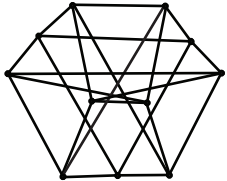}}
\caption{The graph $G_S$ is minor minimal intrinsically knotted. }
\label{GS}
\end{figure}

  It is natural to investigate the graphs obtained by adding one edge to each of the six graphs in ${\mathcal F}\left(K_{7}\right)\setminus{\mathcal F}_\triangle\left(K_{7}\right)$ to see if the additional edge causes the graph to become intrinsically knotted. Goldberg, Mattman, and Naimi \cite{GMN} consider one such graph, denoted $E_9+e$, and 
show that ${\mathcal F}\left(E_9+e\right)$ includes 33 minor minimal intrinsically knotted graphs. In unpublished work, Schwartz~\cite{NS} 
found another minor minimal intrinsically knotted graph, $G_S$, with  22 edges using the same approach (see Figure~\ref{GS}).

In total, there are 92 minor minimal intrinsically knotted graphs with 22 edges: 58 in 
 ${\mathcal F}\left(K_{3,3,1,1}\right)$,  33 in ${\mathcal F}\left(E_9+e\right)$, and
$G_S$.   A computer 
search indicates that there are no other minor minimal intrinsically knotted graphs with 22 edges. In a preprint, Kim et al.~\cite{KLLMO}  verified
this for graphs having either a vertex of valence at least 6 or at least two vertices of valence 5.

Another approach to finding minor minimal intrinsically knotted graphs is to consider graphs with a small number of vertices.  There
are no intrinsically knotted graphs with fewer than 7 vertices, and $K_7$ is the only intrinsically knotted graph with 7 vertices. There are two minor minimal intrinsically knotted graphs with 8 vertices~\cite{BBFFHL}, \cite{CMOPRW}. One is $K_{3,3,1,1}$, the other is the graph obtained by doing a single $\ty$ move on $K_7$.
There are two minor minimal intrinsically knotted graphs with 9 vertices in ${\mathcal F}\left(K_7\right)$, four in ${\mathcal F}\left(K_{3,3,1,1}\right)$, the graph $E_9+e$, and a graph with  28 edges denoted by 
 $G_{9,28}$~\cite{GMN}.  This gives us a total of eight minor minimal intrinsically knotted graphs with 9 vertices. We discuss this result in more detail in a separate article in this volume~\cite{MMR}.

This brings the total number of graphs known to be minor minimal intrinsically knotted to 264. Most appear in the family of
one of the graphs we have mentioned.  In particular, there are
14 in ${\mathcal F}\left(K_7\right)$,
58 in ${\mathcal F}\left(K_{3,3,1,1}\right)$
33  in ${\mathcal F}\left(E_{9}+e\right)$,
and 156 in ${\mathcal F}\left(G_{9,28}\right)$.
The remaining three are $G_S$, $G_{14,25}$~\cite{GMN}, 
and a graph with 13 vertices discovered by Joel Foisy~\cite{F04}.
\medskip

\section{$n$-apex graphs}\label{apex}

In this section we discuss the relationship between a graph being intrinsically linked or knotted and being $n$-apex, which we define below.

\begin{definition} An abstract graph is said to be {\it $n$-apex} if it can be made planar by deleting $n$ or fewer vertices. A $1$-apex graph is also said to be {\it apex}.
\end{definition}

 Robertson and Seymour's Graph Minor Theorem~\cite{rs} implies that for each $n$, there are only finitely many graphs which are minor minimal with respect to being not $n$-apex.  Since any graph which contains a minor that is not $n$-apex, is itself not $n$-apex, a graph is $n$-apex if and only if it contains no graph on that finite list. The graphs $K_5$ and $K_{3,3}$ are the only graphs which are minor minimal with respect to being not $0$-apex.  Sachs \cite{Sa2} proved that all of the graphs in the Petersen family are minor minimal with respect to  being not $1$-apex, and Barsotti and Mattman \cite{BM} proved that these are the only graphs which are not $1$-apex that have 17 or fewer edges.  However, the Petersen family of graphs are not the only graphs which are minor minimal with respect to being not $1$-apex.  We have the following more general results about the negative relationship between $n$-apex and intrinsic linking  and  knotting.

\begin{lemma}
\cite{Sa2}
\label{lemILintrinsically Linked}%
No intrinsically linked graph is $1$-apex.
\end{lemma}

\begin{lemma} \cite{BBFFHL}, \cite{OT}
\label{lemILintrinsically knotted}%
No intrinsically knotted graph is $2$-apex.
\end{lemma}

With this theorem in hand, it makes sense to try to classify all graphs that are minor minimal with respect to being not $2$-apex as a stepping stone to a classification of all minor minimal intrinsically knotted graphs.  However, even the set of graphs which are minor minimal with respect to  being not $1$-apex are not yet characterized.  By doing a computer 
search, Pierce~\cite{P} found a total of 157 minor minimal not $1$-apex graphs.  While this may not be a complete list, it does include all examples that have either 10 or fewer vertices or 21 or fewer edges.  For graphs which are not $2$-apex, Barsotti and Mattman proved the following result.

\begin{theorem} \cite{BM}, \cite{Ma}
\label{thmMMN2A}%
Every graph which is not $2$-apex has at least 21 edges. The 20 graphs in ${\mathcal F} (K_7)$ are the only graphs which are minor minimal with respect to being not $2$-apex and have exactly 21 edges.
\end{theorem}

Another paper in this volume~\cite{MP} describes a computer search for ``small graphs''
 which are minor minimal with respect to being not $2$-apex. In particular, the only such graphs with exactly  22 edges are the 58 graphs in ${\mathcal F}\left(K_{3,3,1,1}\right)$ together with two additional $4$-regular 
graphs with 11 vertices. There are no graphs which are minor minimal with respect to being not $2$-apex that have precisely 23 edges.  On the other hand, for each number between 24 and 30 inclusive, there do exist graphs which are minor minimal with respect to being not $2$-apex which have the specified number of edges (see~\cite{MP} for details).
The same paper shows that there are exactly 12 graphs which are minor minimal with respect to being not $2$-apex and have at most 9 vertices, and all but two of
them (both with 9 vertices) are in ${\mathcal F}\left(K_9\right)\cup{\mathcal F}\left(K_{3,3,1,1}\right)$.

As mentioned above, all of the graphs in the Petersen family are minor minimal with respect to  being not $1$-apex.  Thus all graphs which are minor minimal with respect to  being intrinsically linked are also minor minimal with respect to  being not $1$-apex. By contrast, not every graph which is minor minimal with respect to being intrinsically knotted is minor minimal with respect to being not $2$-apex. For example, none of the 33 minor minimal intrinsically knotted graphs in ${\mathcal F}\left(E_9+e\right)$ is
minor minimal with respect to not being $2$-apex. Of course, by Lemma~\ref{lemILintrinsically knotted}, 
every minor minimal intrinsically knotted graph
must have a minor which is minor minimal with respect to not being $2$-apex. For example, for $E_9+e$,
 the graph
$E_9$ is a proper minor which is minor minimal with respect to not being $2$-apex.   Similarly, 
each of the 33 minor minimal intrinsically knotted graphs in ${\mathcal F}\left(E_9+e\right)$
has a proper minor among the six graphs in ${\mathcal F}\left(K_7\right)\setminus{\mathcal F}_\triangle\left(K_7\right)$ which is minor minimal with respect to not being $2$-apex.

We close this section with a discussion of how not being $1$- or $2$-apex behaves under $\ty$ and $\yt$ moves. First  observe that $K_{3,3}$ is non-planar, yet doing a $\yt$ move on $K_{3,3}$ yields a planar graph.  Thus, the disjoint union of two (resp. three) $K_{3,3}$'s is not $1$-apex (resp. $2$-apex), yet doing a $\yt$ move yields a $1$-apex (resp. $2$-apex) graph. 

On the other hand, consider Jorgensen's graph $J$ (see~\cite{B}) which is minor minimal with respect to not being $1$-apex, but there is a $\ty$ move that renders it $1$-apex.  Indeed, the vertex to be deleted is the one introduced by the $\ty$ move.   Similarly, the disjoint union
$J \sqcup K_5$ is not $2$-apex but becomes $2$-apex with an appropriate $\ty$ move.
In general, the only way a $\ty$ move could fail to preserve the property of
not being $1$-apex or $2$-apex is if the new vertex of valence 3 is the one that is removed (see~\cite{MP}
for more details).

\medskip

\section{Conway-Gordon type theorems for graphs in ${\mathcal F}\left(K_{6}\right)$ and ${\mathcal F}\left(K_{7}\right)$}\label{CGK_6K_7}

In order to prove that $K_6$ is intrinsically linked and $K_7$ is intrinsically knotted, Conway and Gordon obtained the following result about the linking number and Arf invariant of cycles in embeddings of $K_6$ and $K_7$ respectively. Note that ${\rm lk}$ denotes the linking number of a pair of disjoint cycles  and ${\rm Arf}$ denotes the Arf invariant of a cycle.

\begin{ConGor}\label{CG} \cite{CG}
\begin{enumerate}

\item \it{For any spatial embedding $f$ of $K_{6}$, 
$$\sum{\rm lk}\left(f\left(\gamma\right)\right)\equiv 1\pmod{2}$$ where the sum is taken over all pairs of disjoint cycles $\gamma$ in $K_6$. 
\medskip

\item For any spatial embedding $f$ of $K_{7}$, $$\sum{\rm Arf}\left(f\left(\gamma\right)\right)\equiv 1\pmod{2}$$ where  the sum is taken over all $7$-cycles $\gamma$ in $K_7$. }
\end{enumerate}
\end{ConGor}

We would also like to obtain similar results about integer invariants of spatial graphs.  However, we first introduce some notation as follows.  Let $G$ be a graph.  We denote the set cycles in $G$ by $\Gamma\left(G\right)=\Gamma^{(1)}\left(G\right)$, and the set of all $k$-cycles by $\Gamma_{k}\left(G\right)$. We denote the set of all unions of mutually disjoint pairs of cycles of $G$ by $\Gamma^{(2)}\left(G\right)$, and the set of all pairs of a $k$-cycle and a $l$-cycle  by $\Gamma_{k,l}^{(2)}\left(G\right)$. For an element $\gamma$ in $\Gamma^{(r)}\left(G\right)$ and a spatial embedding $f$ of $G$, $f\left(\gamma\right)$ is a (possibly trivial) knot if $r=1$ and a (possibly trivial) $2$-component link if $r = 2$.  A {\it Hamiltonian cycle} of $G$ is a cycle containing every vertex of $G$.  If $\gamma$ is a Hamiltonian cycle, we call $f\left(\gamma\right)$ a {\it Hamiltonian knot} in $f\left(G\right)$ even if $f(\gamma)$ is trivial.   Finally, $a_2(K)$ denotes the second coefficient of the Conway polynomial of $K$.  Note that $a_{2}\left(K\right)\equiv {\rm Arf}\left(K\right) \pmod{2}$.

Nikkuni  \cite{N09b} proved the following theorem about integer invariants of spatial embeddings of $K_6$ and $K_7$.

\begin{theorem}\label{CG_refine} \cite{N09b}
\begin{enumerate}
\item For any spatial embedding $f$ of $K_{6}$, we have
\begin{eqnarray*}
2\sum_{\gamma\in \Gamma_{6}\left(K_{6}\right)}a_{2}\left(f\left(\gamma\right)\right)
-2\sum_{\gamma\in \Gamma_{5}\left(K_{6}\right)}a_{2}\left(f\left(\gamma\right)\right)
=
\sum_{\gamma\in \Gamma^{(2)}\left(K_{6}\right)}{\rm lk}\left(f\left(\gamma\right)\right)^{2}-1 
\end{eqnarray*}. 

\item For any spatial embedding $f$ of $K_{7}$, we have
\begin{eqnarray*}
&&7\sum_{\gamma\in \Gamma_{7}\left(K_{7}\right)}a_{2}\left(f\left(\gamma\right)\right)
-6\sum_{\gamma\in \Gamma_{6}\left(K_{7}\right)}a_{2}\left(f\left(\gamma\right)\right)
-2\sum_{\gamma\in \Gamma_{5}\left(K_{7}\right)}a_{2}\left(f\left(\gamma\right)\right)\\
&=&
2\sum_{\gamma\in \Gamma_{3,4}^{(2)}\left(K_{7}\right)}{\rm lk}\left(f\left(\gamma\right)\right)^{2}-21. 
\end{eqnarray*}
\end{enumerate}
\end{theorem}

 Observe that the Conway-Gordon Theorem can be recovered by taking the mod 2 reduction of Theorem \ref{CG_refine}. Thus we can consider Theorem \ref{CG_refine} to be the integral lift of the Conway-Gordon Theorem.

 Recall from Section~\ref{IKIL} that if a graph $H$ is intrinsically linked (resp. knotted), then any graph $G$ in ${\mathcal F}_{\triangle}\left(H\right)$ is also intrinsically linked (resp. knotted). Moreover, Nikkuni and Taniyama \cite{NT12} established a systematic method of using a Conway-Gordon type equation for $H$ to obtain one for $G$. In particular, by applying their method to the Conway-Gordon Theorem, they obtain the following. 
 

\begin{theorem}\label{NT_main}  \cite{NT12} 
\begin{enumerate}
\item Let $G$ be a graph in ${\mathcal F}_{\triangle}\left(K_{6}\right)$. Then, there exists a map $\omega:\Gamma\left(G\right)\to {\mathbb Z}$ such that for any spatial embedding $f$ of $G$, we have
\begin{eqnarray*}
2\sum_{\gamma\in \Gamma\left(G\right)}\omega\left(\gamma\right)a_{2}\left(f\left(\gamma\right)\right)
=
\sum_{\gamma\in \Gamma^{(2)}\left(G\right)}{\rm lk}\left(f\left(\gamma\right)\right)^{2}
-1.
\end{eqnarray*}

\item Let $G$ be a graph in ${\mathcal F}_{\triangle}\left(K_{7}\right)$. Then, there exists a map $\omega:\Gamma\left(G\right)\cup \Gamma^{(2)}\left(G\right)\to {\mathbb Z}$ such that for any spatial embedding $f$ of $G$, we have
\begin{eqnarray*}
\sum_{\gamma\in \Gamma\left(G\right)}\omega\left(\gamma\right)a_{2}\left(f\left(\gamma\right)\right)
=
2\sum_{\gamma\in \Gamma^{(2)}\left(G\right)}\omega\left(\gamma\right){\rm lk}\left(f\left(\gamma\right)\right)^{2}
-21. 
\end{eqnarray*}
\end{enumerate}
\end{theorem}
\medskip

Hashimoto and Nikkuni  \cite{HN12a} gave specific values for the map $\omega:\Gamma\left(G\right)\to {\mathbb Z}$ in part (1) of Theorem \ref{NT_main}. By taking the mod 2 reduction of Theorem \ref{NT_main}, we have the following. 

\begin{corollary}\label{TY_main_cor}  \cite{Sa2}, \cite{TY01}
\begin{enumerate}
\item Let $G$ be a graph in ${\mathcal F}_{\triangle}\left(K_{6}\right)$. For any spatial embedding $f$ of $G$, we have
$$\sum_{\gamma\in \Gamma^{(2)}\left(G\right)}{\rm lk}\left(f\left(\gamma\right)\right)
\equiv 1 \pmod{2}$$. 

\item Let $G$ be a graph in ${\mathcal F}_{\triangle}\left(K_{7}\right)$. Then, there exists a subset $\Gamma$ of $\Gamma\left(G\right)$ such that for any spatial embedding $f$ of $G$, we have
$$\sum_{\gamma\in \Gamma}{\rm Arf}\left(f\left(\gamma\right)\right) \equiv 1 \pmod{2}$$. 
\end{enumerate}
\end{corollary}

 O'Donnol \cite{D10} showed the following Conway-Gordon type theorem for $K_{3,3,1}$, which is the only graph in ${\mathcal F}\left(K_{6}\right)\setminus {\mathcal F}_{\triangle}\left(K_{6}\right)$. 

\begin{theorem}\label{CG_refine_k331} \cite{D10} For any spatial embedding $f$ of $K_{3,3,1}$, we have
\begin{eqnarray*}
&&2\sum_{\gamma\in \Gamma_{7}\left(K_{3,3,1}\right)}a_{2}\left(f\left(\gamma\right)\right)
-4\sum_{\substack{{\gamma\in \Gamma_{6}\left(K_{3,3,1}\right)} \\ {u\not\in\gamma}}}a_{2}\left(f\left(\gamma\right)\right)
-2\sum_{\gamma\in \Gamma_{5}\left(K_{3,3,1}\right)}a_{2}\left(f\left(\gamma\right)\right)\\
&=&
\sum_{\gamma\in \Gamma_{3,4}^{(2)}\left(K_{3,3,1}\right)}{\rm lk}\left(f\left(\gamma\right)\right)^{2}-1,
\end{eqnarray*}
where $u$ is the unique vertex of $K_{3,3,1}$ with valence $6$. 
\end{theorem}

It now follows that part (1) of Theorem \ref{NT_main} and Corollary \ref{TY_main_cor} hold for any graph in the Petersen family. 

On the other hand, as we have already mentioned the six graphs in ${\mathcal F}\left(K_{7}\right) \setminus {\mathcal F}_{\triangle}\left(K_{7}\right)$ are not intrinsically knotted \cite{FN08}, \cite{HNTY10}, \cite{GMN} (though they are minor-minimal ``intrinsically knotted or completely $3$-linked'', see \cite{HNTY10} for details). Thus for each of these six graphs, there does not exist a subset $\Gamma$ satisfying the condition in part (2) of Corollary \ref{TY_main_cor}.  This leaves us with the following open question.

\begin{question}\label{q1}
Is there an integral Conway-Gordon type formula for every graph in ${\mathcal F}\left(K_{7}\right) \setminus {\mathcal F}_{\triangle}\left(K_{7}\right)$?
\end{question}

\medskip

\section{Conway-Gordon type theorems for $K_{3,3,1,1}$}\label{CGK3311}

 In \cite{MRS88}, Motwani, Raghunathan, and Saran claimed that $K_{3,3,1,1}$ could be shown to be intrinsically knotted  using the same technique as in part (2) of the Conway-Gordon Theorem.  However, Kohara and Suzuki \cite{KS92} showed that there exist two spatial embeddings of $K_{3,3,1,1}$ such that the sum of $a_{2}(\gamma)$ over all of the Hamiltonian knots $\gamma$ for one embedding is even and for the other embedding is odd. Thus the situation for $K_{3,3,1,1}$ is necessarily different from that of $K_{7}$. By using a new approach, Foisy \cite{F02} succeeded in proving that $K_{3,3,1,1}$ is intrinsically knotted. Furthermore, Foisy's proof implies the following. 

\begin{theorem}\label{intrinsically knotted_k3311}\cite{F02}
 For any spatial embedding $f$ of $K_{3,3,1,1}$, there exists a cycle $\gamma$ in $\bigcup_{k=4}^{8}\Gamma_{k}\left(K_{3,3,1,1}\right)$ such that $a_{2}\left(f\left(\gamma\right)\right)\equiv 1\pmod{2}$. 
\end{theorem}

On the other hand, Hashimoto and Nikkuni showed the following Conway-Gordon type theorem for $K_{3,3,1,1}$.  Here, $x$ and $y$ denote the two vertices of $K_{3,3,1,1}$ of valence $7$. 

\begin{theorem}\label{main_theorem_k3311} \cite{HN12b}
\begin{enumerate}
\item For any spatial embedding $f$ of $K_{3,3,1,1}$, we have
\begin{eqnarray*}
&&4\sum_{\gamma\in \Gamma_{8}(K_{3,3,1,1})}a_{2}\left(f\left(\gamma\right)\right)
-4\sum_{\substack{{\gamma\in \Gamma_{7}\left(K_{3,3,1,1}\right)} \\ 
{\left\{x,y\right\}\not\subset \gamma}}}
a_{2}\left(f\left(\gamma\right)\right) \\
&&-4\sum_{\gamma\in \Gamma_{6}'}a_{2}\left(f\left(\gamma\right)\right)
-4\sum_{\substack{{\gamma\in \Gamma_{5}\left(K_{3,3,1,1}\right)} \\ 
{\left\{x,y\right\}\not\subset \gamma}}}
a_{2}\left(f\left(\gamma\right)\right)\\
&=&
\sum_{\gamma\in \Gamma_{3,5}^{(2)}\left(K_{3,3,1,1}\right)}{\rm lk}(f\left(\gamma\right))^{2}
+2\sum_{\gamma\in \Gamma_{4,4}^{(2)}\left(K_{3,3,1,1}\right)}{\rm lk}(f\left(\gamma\right))^{2}
-18, 
\end{eqnarray*}
where $\Gamma'_{6}$ is a specific subset of $\Gamma_{6}\left(K_{3,3,1,1}\right)$ which does not depend on $f$. 
\medskip

\item For any spatial embedding $f$ of $K_{3,3,1,1}$, we have
\begin{eqnarray*}
\sum_{\gamma\in \Gamma_{3,5}^{(2)}\left(K_{3,3,1,1}\right)}{\rm lk}(f\left(\gamma\right))^{2}
+2\sum_{\gamma\in \Gamma_{4,4}^{(2)}\left(K_{3,3,1,1}\right)}{\rm lk}(f\left(\gamma\right))^{2}
\ge 22. 
\end{eqnarray*}
\end{enumerate}
\end{theorem}

By combining parts (1) and (2) of Theorem \ref{main_theorem_k3311}, we have the following Corollary, which refines Theorem \ref{intrinsically knotted_k3311} by identifying the cycles which might be non-trivial knots.

\begin{corollary}\label{main_theorem_k3311_cor} For any spatial embedding $f$ of $K_{3,3,1,1}$, we have \begin{eqnarray}\label{k3311_ineq}
&&\sum_{\gamma\in \Gamma_{8}\left(K_{3,3,1,1}\right)}a_{2}\left(f\left(\gamma\right)\right)
-\sum_{\substack{{\gamma\in \Gamma_{7}\left(K_{3,3,1,1}\right)} \\ 
{\left\{x,y\right\}\not\subset \gamma}}}
a_{2}\left(f\left(\gamma\right)\right) \\
&&-\sum_{\gamma\in \Gamma_{6}'}a_{2}\left(f\left(\gamma\right)\right)
-\sum_{\substack{{\gamma\in \Gamma_{5}\left(K_{3,3,1,1}\right)} \\ 
{\left\{x,y\right\}\not\subset \gamma}}}
a_{2}\left(f\left(\gamma\right)\right)
\ge 1. \nonumber
\end{eqnarray}
\end{corollary}

Kohara and Suzuki's \cite{KS92} two embeddings of $K_{3,3,1,1}$ show that the left side of the above inequality is not necessarily congruent to $1$ modulo $2$. Thus Corollary \ref{main_theorem_k3311_cor} shows that we get useful information by using integer invariants rather than $\mathbb{Z}_2$ invariants. Moreover, the following result of Hashimoto and Nikkuni implies that for any graph $G$ in ${\mathcal F}_{\triangle}\left(K_{3,3,1,1}\right)$, one can obtain an integral Conway-Gordon type inequality from Corollary \ref{main_theorem_k3311_cor} by using Nikkuni and Taniyama's method \cite{NT12}.

\begin{theorem}\label{NT_main_cor} \cite{HN12b} 
Let $G$ be a graph in ${\mathcal F}_{\triangle}\left(K_{3,3,1,1}\right)$. Then, there exists a map $\omega:\Gamma\left(G\right)\to{\mathbb Z}$ such that for any spatial embedding $f$ of $G$, we have
$\sum_{\gamma\in \Gamma\left(G\right)}\omega\left(\gamma\right)a_{2}\left(f\left(\gamma\right)\right)
\ge 1$. 
\end{theorem}

Recall that ${\mathcal F}\left(K_{3,3,1,1}\right)$ contains 58 graphs, 32 of which do not belong to ${\mathcal F}_{\triangle}\left(K_{3,3,1,1}\right)$ but are nonetheless intrinsically knotted \cite{GMN}.  This leaves us with the following open questions.

\begin{question}\label{q2}
Is there an integral Conway-Gordon type formula for every graph in ${\mathcal F}\left(K_{3,3,1,1}\right)\setminus{\mathcal F}_{\triangle}\left(K_{3,3,1,1}\right)$?
\end{question}

\begin{question}\label{q3}
 Are there integral Conway-Gordon type formulae for any of the other minor minimal intrinsically knotted graphs? 
\end{question}

\medskip

\section{Linear embeddings of graphs}\label{linear}

In this section we consider a special type of embedding of graphs in $\mathbb{R}^3$ which is determined entirely by the placement of the vertices.

\begin{definition}A spatial embedding of a graph is said to be {\it linear} (or {\it rectilinear}) if each of the edges in the image is a straight line segment in $\mathbb{R}^3$.  
\end{definition}

The following are noteworthy results about knots and links in linear embeddings of complete graphs.

\begin{theorem}\label{HJ} \cite{Hu},  \cite{HJ07}
Any linear spatial embedding of $K_{6}$ contains at most one trefoil knot and at most three Hopf links. Furthermore,

\begin{enumerate}
\item Such an embedding does not contain a trefoil knot if and only if it contains exactly one Hopf link. 

\item Such an embedding contains a trefoil knot if and only if it contains exactly three Hopf links. 
\end{enumerate}
\end{theorem}

Various authors have obtained results about knots and links in linear embeddings of $K_7$.  In particular, we have.

\begin{theorem}\label{Alfon} 
 \cite{Br},  \cite{RA}
Every linear spatial embedding $K_{7}$ contains a trefoil knot.   
\end{theorem}

\begin{theorem}\label{HuhK7}  \cite{Huh} No linear spatial embedding of $K_7$ contains more than $3$ figure eight knots.
\end{theorem}

\begin{theorem}\label{LinksK7}\cite{Lew}  The number of non-trivial links in any linear spatial embedding of $K_7$ whose convex hull is a polyhedron with seven vertices is between $21$ and $48$.
\end{theorem}

Further results about knots and links in linear embeddings of $K_7$ can be found in the slides from a talk given by Choon Bae Jeon at the International Workshop on Spatial Graphs in 2010 (see \cite{Je}).

\begin{theorem}\label{K9}  \cite{NP}
Every linear spatial embedding of $K_{9}$ contains a non-split 3-component link.   
\end{theorem}

In the 1980's, Sachs \cite{Sa1,Sa2} conjectured that if a graph has an embedding in $\R^3$ with no non-trivial links,
then it has a linear embedding with no non-trivial links.
As far as we know, this conjecture remains open.
By contrast, in spite of Theorem~\ref{K9}, it has been shown that $K_9$ has an embedding with no non-split 3-component links \cite{FNP}. 

Results have also been obtained about knotting and linking in linear embeddings of complete graphs with a large number of vertices.  In particular, we have the following result of Negami.

\begin{theorem}\label{Negami}  \cite{Ne}
For every knot or link $J$, there is an integer $R(J)$ such that every linear embedding of the complete graph $K_{R(J)}$ in $\mathbb{R}^3$ contains $J$.   
\end{theorem}

In order to study the rate of growth of knotting and linking in random linear embeddings of $K_n$ as a function of $n$, we introduce the following definition.

\begin{definition} Let $f:\mathbb{N}\to \mathbb{N}$ be a function of the naturals. Then $f(n)$ is said to be of the {\it order of $\theta(g(n))$}, if there exist constants $c$, $C > 0$ such that for sufficiently large $n$, $$c g(n) \leq f(n) \leq C g(n).$$ \end{definition}

Flapan and Kozai obtained the following results about entanglement in random linear embeddings of $K_n$ in a cube.
 
\begin{theorem}\label{Kenji1} \cite{FK} Let $n\geq 6$, and let $\Gamma$ be the image of a random linear embedding of $K_n$ in a cube.  Then the mean sum of squared linking numbers for pairs of disjoint cycles in $\Gamma$ is of the order of $\theta(n(n!))$.
\end{theorem}
 
\begin{theorem}\label{Kenji2}  \cite{FK}  Let $n\geq 3$, and let $\Gamma$ be a random linear embedding of $K_n$ in a cube.   Then the mean sum of squared writhe for cycles  in $\Gamma$ is of the order of 
	$\theta(n(n!))$.
\end{theorem}

Other significant results that put restrictions on the types and number of knots and links that 
must, can, or cannot occur in linearly embedded graphs include the following result of Naimi and Pavelescu.

\begin{theorem}\label{K331}\cite{NP2} In any linear spatial embedding of $K_{3,3,1}$ containing an odd number of non-trivial links, all such links are Hopf links.  In any linear spatial embedding of $K_{3,3,1}$ containing an even number of non-trivial links, one such link is a $(2,4)$-torus link and the rest are Hopf links.
\end{theorem} 

A variety of approaches have been useful in the study of linear embeddings of specific graphs.  We discuss the use of  Conway-Gordon type theorems and oriented matroid theory in the two subsections that follow.  
\medskip

\subsection{Conway-Gordon type theorems and linear embeddings}

The stick number was introduced to study linear embeddings of knots and links.  In particular, we define the {\it stick number} $s\left(L\right)$ of a link $L$ to be the minimum number of edges in a polygon which represents $L$. Every link contained in the image of a linear spatial embedding of $K_{n}$ has stick number no more than $ n$. 

While the original proofs of Theorems \ref{HJ}, \ref{Alfon}, and \ref{K9} were combinatorial and computational, Nikkuni \cite{N09b} gave much simpler topological proofs of Theorems~\ref{HJ} and \ref{Alfon} by using the {\it stick number} and applying Theorem~\ref{CG_refine}.  In particular, we see as follows that every linear embedding of $K_7$ must contain a trefoil knot.  

Fix a linear spatial embedding $f$ of $K_{7}$.   It is well known that $s\left(K\right)\ge 6$ for any non-trivial knot $K$, moreover, $s\left(K\right)=6$ if and only if $K$ is a trefoil knot, and $s\left(K\right)=7$ if and only if $K$ is a figure eight knot.  Also, $a_{2}\left({\rm trefoil\ knot}\right)=1$ and $a_{2}\left({\rm figure\ eight\ knot}\right)=-1$.  Thus by part (2) of Theorem~\ref{CG_refine}, we have the inequality

\begin{eqnarray*}
7\sum_{\gamma\in \Gamma_{7}\left(K_{7}\right)}a_{2}\left(f\left(\gamma\right)\right)
\ge 
2\sum_{\gamma\in \Gamma_{3,4}^{(2)}\left(K_{7}\right)}{\rm lk}\left(f\left(\gamma\right)\right)^{2}-21. 
\end{eqnarray*}

Fleming and Mellor \cite{FM09} proved that the image of any spatial embedding of $K_{7}$ contains at least $14$ links with odd linking number whose components are a $3$-cycle and a $4$-cycle.  It follows that $$\sum_{\gamma\in \Gamma_{7}\left(K_{7}\right)}a_{2}\left(f\left(\gamma\right)\right) \ge 1.$$ This implies that $f\left(K_{7}\right)$ must contain a trefoil knot. Theorem \ref{HJ} can also be proven using part (1) of Theorem~\ref{CG_refine}, see \cite{N09b}. 

We remark here that it is still unknown whether every spatial embedding of $K_{3,3,1,1}$ must contain a non-trivial Hamiltonian knot or not. But Hashimoto and Nikkuni use Corollary \ref{main_theorem_k3311_cor} to prove that if $f$ is a linear embedding of $K_{3,3,1,1}$, then $f\left(K_{3,3,1,1}\right)$ contains a non-trivial Hamiltonian knot. 

\begin{theorem}\label{lineark3311} \cite{HN12b}
For any linear spatial embedding $f$ of $K_{3,3,1,1}$, there is an $8$-cycle $\gamma$ of $K_{3,3,1,1}$ whose image $f\left(\gamma\right)$ is a non-trivial knot. 
\end{theorem}

We prove Theorem \ref{lineark3311} as follows. Let $G_{x}$ and $G_{y}$ be the subgraphs of $K_{3,3,1,1}$ isomorphic to $K_{3,3,1}$ which are obtained by deleting one of the $7$-valent vertices $x$ and $y$, respectively. Then it follows from Theorem \ref{CG_refine_k331} that $$\sum_{\gamma\in \Gamma_{7}\left(G_{v}\right)}a_{2}\left(f\left(\gamma\right)\right) \ge 0\ (v=x,y).$$ Since the set $\Gamma_{7}\left(G_{x}\right)\cup \Gamma_{7}\left(G_{y}\right)$ consists of $7$-cycles which do not contain both $x$ and $y$, we have the following inequality.

\begin{eqnarray}\label{aaa}
\sum_{\substack{{\gamma\in \Gamma_{7}\left(K_{3,3,1,1}\right)} \\ 
{\left\{x,y\right\}\not\subset \gamma}}}
a_{2}\left(f\left(\gamma\right)\right) 
= \sum_{\gamma\in \Gamma_{7}\left(G_{x}\right)}a_{2}\left(f\left(\gamma\right)\right)
+ \sum_{\gamma\in \Gamma_{7}\left(G_{y}\right)}a_{2}\left(f\left(\gamma\right)\right)
\ge 0. 
\end{eqnarray}

Now by combining (\ref{aaa}) with Corollary \ref{main_theorem_k3311_cor}, we obtain the inequality $$\sum_{\gamma\in \Gamma_{8}\left(K_{3,3,1,1}\right)}a_{2}\left(f\left(\gamma\right)\right)\ge 1.$$ This implies that $f\left(K_{3,3,1,1}\right)$ contains a non-trivial Hamiltonian knot with $a_{2}>0$. We remark here that Jorge Calvo \cite{C01} proved that there are exactly eight knots with $a_{2}>0$ and $s\le 8$.  We close this section with some open questions.

\begin{question}\label{q4}
Does every linear spatial embedding of $K_{3,3,1,1}$ contain a trefoil knot?
\end{question}

\begin{question}\label{q5}
By applying oriented matroid theory (see the next subsection) Naimi and Pavelescu \cite{NP} proved that the number of non-split $2$-component links in any linear spatial embedding of $K_{3,3,1}$ is $1, 2, 3, 4$ or $5$ . Can we give a topological proof by using Theorem~\ref{main_theorem_k3311}? How about for graphs in the Petersen family other than $K_{6}$ and $K_{3,3,1}$ ?
\end{question}

\begin{question}\label{q6}
A graph is said to be {\it intrinsically triple linked} if the image of every spatial embedding of the graph contains a non-split $3$-component link. Although $K_{9}$ is not intrinsically triple linked (Flapan, Naimi, and Pommersheim \cite{FNP}), Naimi and Pavelescu \cite{NP} used oriented matroid theory and a computer to prove Theorem~\ref{K9}. Can we give an alternative proof without the help of a computer by applying a Conway-Gordon type theorem? 
\end{question}
\medskip

\subsection{Oriented matroids and linear embeddings}

Another approach to the study of linear embeddings of graphs is to use oriented matroids.  For the sake of space, we do not give a full definition of an oriented matroid. 
Rather we give enough of a partial definition of
a uniform oriented matroid
 to be able to explain
how they are associated with linear embeddings of graphs.
For a formal description of oriented matroids we refer the reader to \cite{BLSWZB}.

Let $\Gamma$ be a linear spatial graph with $n$ vertices in general position in $\R^3$.  We assign the vertices of $\Gamma$ a fixed (but arbitrary) order.  Now let $a_1$, $a_2$, $a_3$, and $a_4$ be vectors in $\R^3$ which represent four of the vertices of $\Gamma$ such that $a_1<a_2<a_3<a_4$ with respect to our fixed order.  This ordered set of vectors defines a $3 \times 3$ matrix $A = [a_4-a_1 \; | \; a_4-a_2 \; | \; a_4-a_3]$, where each $a_4-a_i$ is a column vector.  Observe that since the vectors $a_1$, $a_2$, $a_3$, $a_4$ are not coplanar, the determinant of $A$ is not zero.  Hence we can assign a $+$ or $-$ sign to the ordered set $\{a_1, a_2, a_3, a_4\}$ according to whether $\det(A)$ is positive or negative.

Each such ordered signed set of four vertices $\{a_1, a_2, a_3, a_4\}$ is said to be a {\it signed basis}, and we define the {\it uniform oriented matroid} $\Omega(\Gamma)$ to be the collection of all such signed bases. Note that $\Omega(\Gamma)$ depends only on the embedding of the vertices of $\Gamma$,
and not in any way on which edges $\Gamma$ contains.  We say that $\Omega(\Gamma)$ is of rank $4$, because we are considering sets of four vertices.  Finally, we let  $\mathrm{OM}(4,n)$ denote the set of all uniform oriented matroids of rank $4$ of linear spatial graphs with $n$ vertices.  

We will see below that one can compute
the linking number of any pair of disjoint cycles in $\Gamma$ 
directly from  $\Omega(\Gamma)$,
without specifying the embedding $\Gamma$ itself.
Thus, one could study  ``linking behavior'' in all  linear embeddings of 
a graph with $n$ vertices
 if one had a list of all oriented matroids in $\mathrm{OM}(4,n)$ ---
without  having a list of all the linear embeddings of the graph.  For example, in Theorem~\ref{K9} it was determined that
every linear embedding of $K_9$ contains a non-split 3-component link.
This was done by going through all of the elements of $\mathrm{OM}(4,9)$
and checking that there is always a 3-component link such that 
one of its components has non-zero linking number 
with each of the other two.

Catalogs of uniform oriented matroids for various ranks and numbers of vertices have been found using a computer, 
and are available on the ``Homepage of Oriented Matroids" \cite{Fi}.
However, as of this writing, for rank 4 the list only goes up to $n = 9$ vertices.

\begin{figure}[h]
{\includegraphics{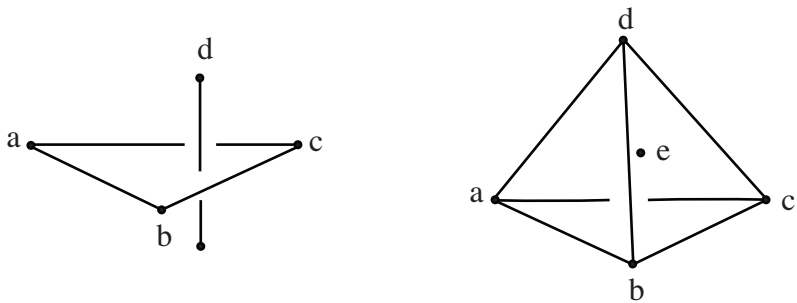}}
\caption{
On the left, the edge $\{d,e\}$ pierces the triangle $\{a,b,c\}$. 
On the right, vertex $e$ lies inside the tetrahedron $\{a,b,c,d\}$.  }
\label{circuits}
\end{figure}

A different but equivalent way to associate an oriented matroid with $\Gamma$ 
is via circuits instead of bases.
For each collection of 5 of the $n$ vertices of $\Gamma$,
either one of the vertices will be inside the tetrahedron determined by the remaining four vertices,
or two of the vertices will determine an edge that intersects 
the interior of the disk bounded by the triangle determined by the remaining three vertices.  The two possibilities are illustrated in Figure~\ref{circuits}. 
According to which of these possibilities occurs for a given subset of 5 vertices, the subset  is assigned a 4--1 or a 3--2 partition.
These partitions are then sufficient to describe an oriented matroid of rank 4 on $n$ elements.  More generally, a uniform oriented matroid of rank $r$ on $n$ vertices can be defined by assigning a partition to each of the $(r+1)$-subsets (called circuits) of its elements, subject to certain conditions (called circuit axioms).


The oriented matroids cataloged in \cite{Fi} are each given as a list of signed bases,
rather than a list of partitions.
We can obtain the partitions from the signed bases (and vice versa) 
using a simple procedure.
Once we obtain partitions of all subsets of five vertices, we will know which edges pierce which triangles.  
This is sufficient for computing the linking number between any pair of disjoint triangles.
From this, we can compute the linking number between any  two disjoint cycles 
by writing each cycle as a ``sum'' of triangles 
(regardless of whether or not each of the triangles is a subgraph of $\Gamma$).
A detailed example of how to obtain 3--2 partitions from signed bases
and compute linking numbers 
can be found in \cite{NP-manual}.

\begin{figure}[h]
 \includegraphics{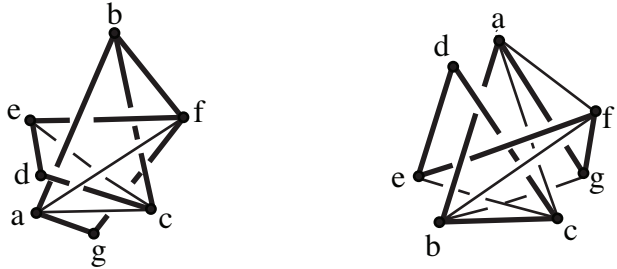}
 \caption{Two linear embeddings of $K_7$, each with a trefoil knot highlighted.}
 \label{figfromRA}
\end{figure}

As an example,
we give a brief outline of the proof of Theorem~\ref{Alfon} by Ramirez Alfonsin \cite{RA} which used oriented matroids to show that 
every linear embedding of $K_7$ contains a trefoil.
Let $M$ be an oriented matroid with vertices, $a$, $b$, $c$, $d$, $e$, $f$, $g$,
that contains the 3--2 partitions 
$(abc,ef)$, $(afb,cd)$, $(cde, ab)$,  $(bcd, ef)$;
and does not contain the 3-2 partitions
$(bcd,ag)$, $(bcd ,fg)$, $(def, ag)$, $(cde, ag)$, $(cde, fg)$.
This restricts the seven vertices to be in one of the two configurations
depicted in Figure~\ref{figfromRA}.

As shown in the figure, 
each of the configurations contains a trefoil.
Now, in addition to the condition on 3--2 partitions given above, 
two other similar conditions are listed in \cite{RA}
which also each guarantee the existence of a trefoil.
Then, using a computer program, it was verified
that every oriented matroid in $\mathrm{OM}(4,7)$
satisfies at least one of the three conditions.

\medskip

\section{Symmetries of spatial graphs in $S^3$}\label{symmetries}

As mentioned at the beginning of this survey, the study of symmetries of spatial graphs was originally motivated by the need to describe the symmetries of non-rigid molecules.  However, just as the characterization of the symmetries of knots and links plays an important role in knot theory, the characterization of symmetries of spatial graphs helps us distinguish  different embeddings of a given graph. 

\begin{definition} Let $\gamma$ be an abstract graph with automorphism group $\Aut(\gamma)$, and let $\Gamma$ be the image of an embedding of $\gamma$ in $S^3$.  The {\em topological symmetry group} of $\Gamma$, denoted by $\mathrm{TSG}(\Gamma)$, is the subgroup of $\Aut(\gamma)$ induced by homeomorphisms of the pair $(S^3, \Gamma)$.  If we restrict consideration to orientation preserving homeomorphisms of the pair $(S^3,\Gamma)$, then we obtain $\TSG(\Gamma)$. \end{definition}

 In this survey we will restrict our attention to $\TSG(\Gamma)$.  However, we abuse notation and will refer to this group   as the {\em topological symmetry group} of $\Gamma$.

Fruct \cite{fr} showed that every finite group can be realized as $\Aut(\gamma)$ for some graph $\gamma$. By contrast, Flapan, Naimi, Pommersheim, and Tamvakis \cite{fnpt} proved that no alternating group $A_n$ with $n>5$ can be the topological symmetry group of a spatial graph in $S^3$.  Furthermore, for $3$-connected graphs, they prove the following.  Note that a graph is {\it $3$-connected} if at least three vertices have to be removed in order to disconnect it or reduce it to a single vertex.

\begin{theorem} \label{isometry}\cite{fnpt}
\label{SO(4)} 
The topological symmetry group of every $3$-connected spatial graph is isomorphic to a finite subgroup of $\mathrm{SO}(4)$.
\end{theorem}

The finite subgroups of $\so(4)$ have been classified and can all be described as quotients of products of cyclic groups $\Z_m$, dihedral groups $D_m$, and the symmetry groups of the regular polyhedra ($A_4$, $S_4$ and $A_5$) \cite{du}. However, Theorem \ref{isometry} does not give any information as to which graphs can be used to realize which groups.  Subsequent research on topological symmetry groups of spatial graphs has focused on particular families of graphs, and determined which graphs in the family can be used to realize particular groups.
\medskip

\subsection{Complete Graphs $K_n$}

The first family of graphs whose topological symmetry groups have been studied are the complete graphs.  Flapan, Naimi and Tamvakis proved the following theorem classifying all groups which can be realized as the topological symmetry group of an embedding of a complete graph.

\begin{theorem} \cite{fnt}
\label{T:TSG2} 
A finite group $H$ is isomorphic to $\TSG(\Gamma)$ for the image of some embedding $\Gamma$ of a complete graph in $S^3$ if and only if $H$ is a finite cyclic group, a dihedral group, a subgroup of $D_m \times D_m$ for some odd $m$, or $A_4$, $S_4$, or $A_5$.
\end{theorem}

Subsequently, in a series of papers, Flapan, Mellor, Naimi, and Yoshizawa determined exactly which complete graphs have embeddings that realize each of the above groups.  In particular, we have the following results.

\begin{theorem} \cite{fmn2}
A complete graph $K_n$ with $n\geq 4$ has an embedding in $S^3$ with image $\Gamma$ such that  \begin{itemize}
	\item $\TSG(\Gamma) = A_4$ if and only if $n \equiv 0$, $1$, $4$, $5$, $8 \pmod {12}$.
	\item $\TSG(\Gamma) = A_5$ if and only if $n \equiv 0$, $1$, $5$, $20 \pmod{60}$.
	\item $\TSG(\Gamma) = S_4$ if and only if $n \equiv 0$, $4$, $8$, $12$, $20 \pmod {24}$.
\end{itemize}
\end{theorem}

\begin{theorem} \label{T:Dm} \cite{fmn3}
A complete graph $K_n$ with $n > 6$ has an embedding in $S^3$ with image $\Gamma$ such that $\TSG(\Gamma) =G$ where $G= \Z_m$ or $D_m$  if and only if one of the following conditions holds:
\begin{enumerate}
	\item $m \geq 4$ is even, and $n \equiv 0 \pmod m$.
	\item $m \geq 3$ is odd and $n \equiv 0, 1, 2, 3 \pmod m$.
	\item $G= D_2$, and $n \equiv 0, 1, 2 \pmod 4$.
	\item $G = \Z_2$.
\end{enumerate}
\end{theorem}

\begin{theorem} \label{T:ZxZxZ} \cite{fmn3}
A complete graph $K_n$ with $n > 6$ has an embedding in $S^3$ whose image  $\Gamma$ has $\TSG(\Gamma) = G$ where $G =  \Z_r \x \Z_s$ or $(\Z_r \x \Z_s) \rtimes \Z_2$ where $r, s$ are odd and $\gcd(r,s) > 1$ if and only if one of the following conditions holds:
\begin{enumerate}
	\item $rs \vert n$.
	\item $\gcd(r,s) = 3$ and $rs \vert (n-3)$.
	\item $G = \Z_3 \x \Z_3$ and $9 \vert (n-6)$.
	\item $G = (\Z_3 \x \Z_3) \rtimes \Z_2$ and $18 \vert (n-6)$.
\end{enumerate}
\end{theorem}

\begin{theorem} \label{T:DxD} \cite{fmn3}
A complete graph $K_n$ with $n > 6$ has an embedding $S^3$ whose image $\Gamma$ has $\TSG(\Gamma) = G$ where $G =  \Z_r \x D_s$ or $D_r \times D_s$ and $r, s \geq 3$ are odd if and only if one of the following conditions holds:
\begin{enumerate}
	\item $2rs \vert n$.
	\item $G = \Z_3 \x D_3$ and $18 \vert (n-6)$.
	\item $G = D_3 \times D_3$ and $36 \vert (n-6)$.
\end{enumerate}
\end{theorem}

The last three results relied on Flapan's earlier classification \cite{fl} of which automorphisms of a complete graph $K_n$ with $n>6$ could be realized by a homeomorphism of $(S^3, \Gamma)$ for some embedding of $K_n$ with image $\Gamma$.  For $n\leq 6$, the automorphism $(1234)$ of $K_6$ is the only automorphism of $K_n$ that cannot be realized by a homeomorphism of $(S^3, \Gamma)$ for some embedding of $K_n$ with image $\Gamma$.  Chambers and Flapan \cite{cf} completed the classification of topological symmetry groups of compete graphs by proving the theorem below which determines all groups that can be realized as the topological symmetry group of some embedding of $K_n$ for $n \leq 6$.  Note for the groups $S_4$, $A_4$, and $A_5$ the classification was already known from \cite{fmn2}.

\begin{theorem} \label{T:nleq6} \cite{cf}
Let $3\leq n \leq 6$.  A non-trivial group $G$ can be realized as $\TSG(\Gamma)$ for some embedding of $K_n$ in $S^3$ with image $\Gamma$ if and only if one of the following conditions holds:
\begin{enumerate}
	\item $n=3$ and $G=\Z_3$ or $D_3$.
	\item $n=4$ and $G=S_4$, $A_4$, $\Z_m$, or $D_m$ where $2$, $3$, or $4$.
	\item $n=5$ and $G=A_5$, $A_4$, $\Z_m$, or $D_m$ where $m=2$, $3$, or $5$.
	\item $n=6$ and $G=D_3\times D_3$, $D_3\times \Z_3$, $\Z_3\times \Z_3$, $(\Z_3\times \Z_3)\x\Z_2$ $\Z_m$, or $D_m$ where $m=2$, $3$, $5$ or $6$.
	 \end{enumerate}
\end{theorem}

\medskip

\subsection{Complete Bipartite Graphs $K_{n,n}$}

Unlike the complete graphs, where only some of the subgroups of $\so(4)$ are realizable as topological symmetry groups, {\em any} finite subgroup of $\so(4)$ can be realized as the topological symmetry group of an embedding of some $K_{n,n}$ \cite{fnpt}.  So the complete bipartite graphs are a natural family of graphs to investigate in order to  understand topological symmetry groups more generally.  

Flapan, Lehle, Mellor, Pittluck and Vongsathorn \cite{flmpv} took the first step in this investigation by classifying which automorphisms of a complete bipartite graph $K_{n,n}$ can be realized by a homeomorphism of $(S^3, \Gamma)$ where $\Gamma$ is the image of some embedding of $K_{n,n}$ in $S^3$.  Mellor proved the following theorem classifying which complete bipartite graphs have embeddings whose topological symmetry groups are  $A_4$, $S_4$ or $A_5$.  

\begin{theorem} \cite{me}
A complete bipartite graph $K_{n,n}$ has an embedding in $S^3$ with image $\Gamma$ such that    \begin{itemize}
	\item $\TSG(\Gamma) = A_4$ if and only if $n \equiv 0$, $2$, $4$, $6$, $8 \pmod {12}$ and $n \geq 4$.
	\item $\TSG(\Gamma) = S_4$ if and only if $n \equiv 0$, $2$, $4$, $6$, $8 \pmod {12}$, $n \geq 4$ and $n \neq 6$.
	\item $\TSG(\Gamma) = A_5$ if and only if $n \equiv 0$, $2$, $12$, $20$, $30$, $32$, $42$, $50 \pmod{60}$ and $n > 30$.
\end{itemize}
\end{theorem}

Hake, Mellor and Pittluck \cite{hmp} considered the subgroups $\Z_m$, $D_m$, $\Z_r \times \Z_s$ and $(\Z_r \times \Z_s) \ltimes \Z_2$ of $\so(4)$.  Their results are summarized in the following theorems:

\begin{theorem}\label{T:cyclic} \cite{hmp} A complete bipartite graph $K_{n,n}$ with $n > 2$ has an embedding in $S^3$ whose image  $\Gamma$ has $\TSG(\Gamma)=G$  for $G=\mathbb{Z}_{m}$ or  $D_{m}$ if and only if one of the following conditions hold:
	\begin{enumerate}
	\item $n\equiv 0,1,2 \pmod{m}$,
	\item $n\equiv 0 \pmod{\frac{m}{2}}$ when $m$ is even,
	\item $n\equiv 2 \pmod{\frac{m}{2}}$ when $m$ is even and $4|m$.
\end{enumerate}
\end{theorem}

\begin{theorem}\label{T:product} \cite{hmp}
A complete bipartite graph $K_{n,n}$ with $n > 2$ has an embedding in $S^3$ whose image  $\Gamma$ has $G\leq \TSG(\Gamma)$ for $G = \mathbb{Z}_{r}\times\mathbb{Z}_{s}$ or $(\mathbb{Z}_{r}\times\mathbb{Z}_{s})\ltimes\mathbb{Z}_{2}$, where $r \vert s$, if and only if one of the following conditions hold:
	\begin{enumerate}	
	\item $n\equiv 0 \pmod {s}$,
	\item $n\equiv 2 \pmod {2s}$ when $r=2$,
	\item $n\equiv s+2 \pmod {2s}$ when $4\vert s$, and $r=2$,
	\item $n\equiv 2 \pmod {2s}$ when $r=4$.
	\end{enumerate}
\end{theorem}

In fact, in each of the cases in Theorem~\ref{T:product}, we can choose the embedding such that $\TSG(\Gamma) = G$ except in the following cases, which are still open: \begin{itemize}
	\item $K_{ls, ls}$, when $1 \leq l < 2r$, $G = \mathbb{Z}_{r}\times\mathbb{Z}_{s}$ or $(\mathbb{Z}_{r}\times\mathbb{Z}_{s})\ltimes\mathbb{Z}_{2}$
	\item $K_{6,6}$, when $G = (\Z_2 \x \Z_4) \ltimes \Z_2$
	\item $K_{10,10}$, when $G = (\Z_4 \x \Z_4) \ltimes \Z_2$
\end{itemize}

The classification of topological symmetry groups for $K_{n,n}$ is far from complete.  In particular, there are many other subgroups of $\so(4)$ to be considered.
\medskip

\subsection{M\"{o}bius ladders}

The other family of graphs whose topological symmetry groups have been classified are the {\it M\"{o}bius ladders} denoted by $M_n$, consisting of a $2n$-cycle together with edges joining each pair of antipodal vertices.  These graphs are referred to as  M\"{o}bius ladders because they can be embedded in $S^3$ in the form of a ladder with $n$ rungs whose ends are joined with a half-twist to resemble a M\"{o}bius strip.  Flapan and Lawrence classified all topological symmetry groups of embeddings of M\"{o}bius ladders in the following two theorems.

\begin{theorem} \cite{fla} The M\"{o}bius ladder $M_3 = K_{3,3}$ has an embedding in $S^3$ whose image  $\Gamma$ has $G\leq \TSG(\Gamma)$ if and only if $G$ is  $D_6$, $D_3$, $D_2$, $\Z_6$, $\Z_3$, $\Z_2$, $D_3 \x D_3$, $\Z_3 \x \Z_3$, $(\Z_3 \x \Z_3) \ltimes \Z_2$, or $D_3 \x \Z_3$.
\end{theorem}

\begin{theorem} \cite{fla} A M\"{o}bius ladder $M_n$ with $n > 3$ has an embedding in $S^3$ whose image  $\Gamma$ has $G\leq \TSG(\Gamma)$ if and only if $G$ is a subgroup of $D_{2n}$.
\end{theorem}
\medskip

\subsection{The mapping class group}

We see as follows that the topological symmetry group is not the only group which represents the symmetries of a spatial graph.  In particular, we can define the mapping class group of a spatial graph as follows.

\begin{definition} Let $\Gamma$ be the image of an embedding of a graph in $S^3$.  Then the mapping class group $\mathrm{MCG}(S^3,\Gamma)$  is the group of isotopy classes of homeomorphisms of the pair $(S^3,\Gamma)$.  If we restrict consideration to orientation preserving homeomorphisms of the pair $(S^3,\Gamma)$, then we obtain $\mathrm{MCG}_+(S^3,\Gamma)$ 
\end{definition}

While the definition of the mapping class group is quite different from that of the topological symmetry group, Cho and Koda \cite{CK} proved that in many cases these groups are isomorphic.

\begin{theorem} \label{MCG}  \cite{CK} Let $\Gamma$ be the image of an embedding of a $3$-connected graph in $S^3$.    Then $\mathrm{MCG}_+(S^3,\Gamma)\cong \mathrm{TSG}_+(S^3,\Gamma)$ if and only if the complement of a neighborhood of $\Gamma$ is atoroidal.  
\end{theorem}

\begin{theorem} \label{MCG}  \cite{CK} Let $\Gamma$ be the image of an embedding of a graph in $S^3$ which is not a knot, and suppose that $\mathrm{MCG}_+(S^3,\Gamma)\cong \mathrm{TSG}_+(S^3,\Gamma)$.  Then $\mathrm{MCG}(S^3,\Gamma)\cong \mathrm{TSG}(S^3,\Gamma)$ if and only if $\Gamma$ is not a planar embedding.
\end{theorem}

\medskip

\subsection{Achirality}
The last type of symmetry that we consider is mirror image symmetry of spatial graphs.  In particular, we have the following definition.

\begin{definition} A spatial graph $\Gamma$ in $S^3$ is said to be {\it chiral} if there does not exist an orientation reversing homeomorphism of $S^3$ which takes $\Gamma$ to itself.  The spatial graph $\Gamma$ is said to be {\it achiral} if such a homeomorphism does exist.  
\end{definition}

While the chirality of knots and links has been studied extensively; relatively little work has been done on the chirality of spatial graphs.  Any graph containing a cycle has a chiral embedding obtained by tying chiral knots in selected edges.  The question is whether a given graph has any embedding which is achiral. If no such embedding exists we say the graph is {\it intrinsically chiral}.  

Flapan \cite{fl1} gave the first examples of intrinsically chiral graphs by showing that a M\"{o}bius ladder with $n \geq 3$ rungs is intrinsically chiral if and only if $n$ is odd.  Flapan and Weaver \cite{fw} used this result to help show that a complete graph $K_n$ is intrinsically chiral if and only if $n = 4k+3$ for some $k \geq 1$.  Flapan, Fletcher, and Nikkuni \cite{ffn} gave an alternative proof of these results by using generalized Simon invariants.  They were also able to use this method to prove that the Heawood graph is intrinsically chiral.  Also, Flapan and Fletcher \cite{ff} classified which complete multipartite graphs admit achiral embeddings.  Finally, Flapan \cite{fl} showed that any non-planar graph which does not have an order 2 automorphism is intrinsically chiral.

The above examples illustrate that, in contrast with the properties of intrinsic knotting or linking which are inherited from minors, a graph which has an intrinsically chiral minor is not in general intrinsically chiral.  There is no known classification of all intrinsically chiral graphs, and the fact that intrinsic chirality is not inherited from minors makes it seem like a classification will be quite difficult.


\medskip

\section{Graphs embedded in $3$-Manifolds}\label{3man}

In this final section, we present some results about spatial graphs embedded in $3$-manifolds.  We would like to study intrinsically knotted and linked graphs for spatial graphs in arbitrary $3$-manifolds.  However, first we need to generalize the concept of a knot and a link.  We use the following definition, though this is certainly not the only way to define a trivial knot and link in a $3$-manifold (see for example \cite{Bus}). 
 
 \begin{definition} Let $K$ and $J$ be disjoint simple closed curves in a $3$-manifold $M$.  We say $K$ is a {\it trivial knot} if $K$ bounds a disk in $M$.  We say $K\cup J$ is a {\it trivial link} if $K$ and $J$ bound disjoint disks in $M$ \end{definition}
 
 Using this definition, Flapan, Howards, Lawrence, Mellor \cite{FHLM} obtained the following somewhat surprising result. 
 
 \begin{theorem}  \cite{FHLM}  A graph is intrinsically linked (resp. intrinsically knotted) in $S^3$ if and only if the graph is intrinsically linked (resp. intrinsically knotted) in any $3$-manifold.
 \end{theorem}
 
 Note that the proof of the intrinsic knotting part of this result relies on the Poincar\'{e} Conjecture \cite{mf}.  
 
 In contrast with the above theorem, the symmetries of a spatial graph in $S^3$ are not generally the same as its symmetries in other $3$-manifolds.  For example, graphs which are intrinsically chiral in $S^3$ are not necessarily intrinsically chiral in other $3$-manifolds. In particular, Flapan and Howards prove the following result.

 \begin{theorem} \cite{fh} Every graph has an achiral embedding in infinitely many closed, connected, orientable, irreducible 3-manifolds. 
 \end{theorem}
 
    On the other hand, they also prove the following.
    
    \begin{theorem}\cite{fh}
For any closed, connected, orientable, irreducible 3-manifold $M$, there are infinitely many graphs which are intrinsically chiral in $M$.
 \end{theorem}

 We can also generalize the definition of the topological symmetry group of a spatial graph in $S^3$ to a spatial graph in any $3$-manifold.  
 
 \begin{definition}  If $\Gamma$ is the image of an embedding of a graph $\gamma$ in a $3$-manifold $M$, we define $\mathrm{TSG}(M,\Gamma)$ to be the subgroup of $\mathrm{Aut}(\gamma)$ induced by homeomorphisms of the pair $(M,\Gamma)$.  By restricting to orientation preserving homeomorphisms of $M$, we obtain $\mathrm{TSG}_{+}(M,\Gamma)$.
 \end{definition}
 
 \begin{theorem}  \cite{ft} Let $M$ be a closed, connected, orientable, irreducible 3-manifold. Then there is an alternating group $A_n$ which is not isomorphic to $\mathrm{TSG}(M,\Gamma)$ for any spatial graph $\Gamma$  in $M$.
 \end{theorem}
 
 Recall from Section~\ref{symmetries} that for any $n>5$, the alternating group $A_n$  is not isomorphic to $\mathrm{TSG}(S^3,\Gamma)$ for any spatial graph $\Gamma$ in $S^3$. The above theorem is for a fixed 3-manifold $M$.  But if we allow $M$ to vary, then every finite group can occur.  In particular, Flapan and Tamvakis prove the following theorem.
 
 \begin{theorem} \cite{ft} For every finite group $G$, there is a hyperbolic rational homology sphere $M$ and a 3-connected spatial graph $\Gamma$  in $M$ such that $G=\mathrm{TSG}(M,\Gamma)$.
 \end{theorem}
 
 We saw in the last section that for any $3$-connected spatial graph $\Gamma$ in $S^3$, $\mathrm{TSG}(S^3,\Gamma)$ is isomorphic to a finite subgroup of $\so(4)$.  By contrast, we have the following result for spaces which are not Seifert fibered.
 
 \begin{theorem} \cite{ft}
 For every closed, orientable, irreducible, $3$-manifold $M$ which is not Seifert fibered, there is a $3$-connected spatial graph $\Gamma$ in $M$ such that $\mathrm{TSG}_+(M,\Gamma)$ is not isomorphic to any group of orientation preserving homeomorphisms of $M$.
\end{theorem}

 \medskip
 
     For more results about symmetries of spatial graphs in $S^3$ and other $3$-manifolds see \cite{Ik1}, \cite{Ik2}, \cite{Ik3}, \cite{Ik4}, \cite{Ko}, \cite{NT09}.

\small
\bibliographystyle{amsplain}

\normalsize

\end{document}